\documentclass{amsart}
\usepackage{amsfonts,amscd,amssymb,amsmath,amsthm}
\usepackage{graphicx,mathrsfs,appendix}
\usepackage{centernot,caption,subcaption,color,url}
\usepackage[T1]{fontenc}
\usepackage{ctable}
\begin{document}

\newtheorem{theorem}{Theorem}[section]
\newtheorem{lemma}[theorem]{Lemma}
\newtheorem{corollary}[theorem]{Corollary}
\newtheorem{conjecture}[theorem]{Conjecture}
\newtheorem{cor}[theorem]{Corollary}
\newtheorem{proposition}[theorem]{Proposition}
\theoremstyle{definition}
\newtheorem{definition}[theorem]{Definition}
\newtheorem{example}[theorem]{Example}
\newtheorem{claim}[theorem]{Claim}
\newtheorem{remark}[theorem]{Remark}

\newenvironment{pfofthm}[1]
{\par\vskip2\parsep\noindent{\sc Proof of\ #1. }}{{\hfill
$\Box$}
\par\vskip2\parsep}
\newenvironment{pfoflem}[1]
{\par\vskip2\parsep\noindent{\sc Proof of Lemma\ #1. }}{{\hfill
$\Box$}
\par\vskip2\parsep}


\newcommand{\R}{\mathbb{R}}
\newcommand{\T}{\mathcal{T}}
\newcommand{\C}{\mathcal{C}}
\newcommand{\B}{\mathcal{B}}
\newcommand{\A}{\mathcal{A}}
\newcommand{\G}{\mathcal{G}}
\newcommand{\Z}{\mathbb{Z}}
\newcommand{\Q}{\mathbb{Q}}
\newcommand{\E}{\mathbb E}
\newcommand{\N}{\mathbb N}

\newcommand{\barray}{\begin{eqnarray*}}
\newcommand{\earray}{\end{eqnarray*}}

\newcommand{\beq}{\begin{equation}}
\newcommand{\eeq}{\end{equation}}


\renewcommand{\Pr}{\mathbb{P}}
\newcommand{\as}{\text{a.s.}}
\newcommand{\Prob}{\Pr}
\newcommand{\Exp}{\mathbb{E}}
\newcommand{\expect}{\Exp}
\newcommand{\1}{\mathbf{1}}
\newcommand{\prob}{\Pr}
\newcommand{\pr}{\Pr}
\newcommand{\filt}{\mathcal{F}}
\newcommand{\F}{\mathcal{F}}
\newcommand{\Bernoulli}{\operatorname{Bernoulli}}
\newcommand{\Binomial}{\operatorname{Binom}}
\newcommand{\Beta}{\operatorname{Beta}}
\newcommand{\Binom}{\Binomial}
\newcommand{\Poisson}{\operatorname{Poisson}}
\newcommand{\Exponential}{\operatorname{Exp}}


\newcommand{\link}{\mbox{lk}}
\newcommand{\Deg}{\operatorname{deg}}
\newcommand{\vertexsetof}[1]{V\left({#1}\right)}
\renewcommand{\deg}{\Deg}
\newcommand{\oneE}[2]{\mathbf{1}_{#1 \leftrightarrow #2}}
\newcommand{\ebetween}[2]{{#1} \leftrightarrow {#2}}
\newcommand{\noebetween}[2]{{#1} \centernot{\leftrightarrow} {#2}}
\newcommand{\Gap}{\ensuremath{\tilde \lambda_2 \vee |\tilde \lambda_n|}}
\newcommand{\dset}[2]{\ensuremath{ e({#1},{#2})}}
\newcommand{\EL}{{ L}}
\newcommand{\ER}{{Erd\H{o}s--R\'{e}nyi }}
\newcommand{\zuk}{{\.{Z}uk}}


\newcommand{\frm}{\ensuremath{ 2\log\log m}}
\newcommand{\csubzero}{c_{0}}
\newcommand{\csubone}{c_{1}}
\newcommand{\csubtwo}{c_{2}}
\newcommand{\csubthree}{c_{3}}
\newcommand{\csubstar}{c_{*}}
\newcommand{\INE}{I^{\epsilon}}
\newcommand{\rsp}{1-C\exp(-md^{1/4}\log n)}
\newcommand{\lc}{\ensuremath{ \operatorname{light}(x,y)}}
\newcommand{\hc}{\ensuremath{ \operatorname{heavy}(x,y)}}

\title[Sublinear preferential attachment with choice]{Sublinear preferential attachment combined with the growing number of choices}
\author{Yury Malyshkin}
\address{Tver State University\\ Moscow Institute of Physics and Technology}
\email{yury.malyshkin@mail.ru}
\date{\today}

\begin{abstract}
We prove almost sure convergence of the maximum degree in an evolving graph model combining a growing number of local choices with sublinear preferential attachment. At each step in the growth of the graph, a new vertex is introduced. Then we draw a random number of edges from it to existing vertices, chosen independently by the following rule. For each edge, we consider a sample of the growing size of vertices chosen with probabilities proportional to the sublinear function of their degrees. Then new vertex attaches to the vertex with the highest degree from the sample. Depending on the growth rate of the sample and the sublinear function, the maximum degree could be of the sublinear order, of the linear order or having almost all edges drawing to it. The prove using various stochastic approximation processes and a large deviation approach.   
\end{abstract}

\maketitle

\section{Introduction}

Preferential attachment graphs are used to model different complex network that exibited certain properties, in particular power law degree distribution. The standart preferential attachment graph, introduced in \cite{barabasi} by Barab{\'a}si and Albert, is constracted by following way. We start with some initial graph $G_0$ on $n_0$ vertices $v_{1-n_0},...,v_0$. Then, the graph $G_{n+1}$ is build from $G_n$ by adding new vertex $v_{n+1}$ and drawing $m$ edges from it to already existing vertices $Y^n_1,...,Y^n_m\in\{v_{1-n_0},...,v_n\}$ chosen independently from each other with probabilities proportion to their degrees, i.e.
$$\Pr(Y^n_i=v_j)=\frac{\deg_{G_n} v_j}{\sum_{k=1-n_0}^{n}\deg_{G_n} v_k}.$$
For this model, many of it properties have been studied. In present paper we are interested in degree distribution and maximum degree of the modification of this model. 
Since asymptotic degree distribution of preferential attachment graph does not depend on initial graph, for simplification of the formulas it is usually suggested that we start with graph, that consist of the single vertex, i.e. $G_0$ consists of vertex $v_0$ and $G_1$ consists of vertices $v_0,v_1$ and $m$ edges between them. 

There are different way to generalize and modify standart preferential attachment model. One of them is to use increasing weighted function $w(x)$, so the vertices chosen with probability proportional to the function of the degree:
$$\Pr(Y^n_i=v_j)=\frac{w(\deg_{G_n} v_j)}{\sum_{k=1-n_0}^{n}w(\deg_{G_n} v_k)}.$$
The linear case was studied in \cite{Mori02,Mori05}, where M{\'o}ri proved that for $w(x)=x+c$, $c>-1$ and $m=1$ ($m=1$ was considered for simplification) the degree distribution follow power law with power $-(3+c)$ and maximum degree is of order $n^{\frac{1}{2+c}}$.
The case of nonlinear weighted function of form $w(x)=x^{\alpha}$ (with $\alpha>0$) was studied in \cite{Athreya08}. For sublenear case ($\alpha<1$) the degree distribution has exponential tails and maximum degree is of order $(\ln n)^{b}$ for some $b>0$ and for superlinear case ($\alpha>1$) the degree distribution is degenerate and maximum degree is asimptotically $n$.

The other way to generalize the model is the addition of choice. In this case, when new vertex added to the graph, it first selects a sample of vertices and then attaches to one of them according to some rule. There were considered different types of this rule, for example in \cite{KR14, MP15, HJ16} authors used rule based on the degree of the vertices and in \cite{HJY18} location-based choice have been used. The effect of the choice (from the sample of $d$ independently chosen vertices) is somewhat similar to the effect of the nonlinear weighted function. In case of min choice, as was shown in \cite{MP15}, maximum degree asymptotically $\ln\ln n/\ln d$ and for max choice and linear weighted function maximum degree could be made both of sublinear and linear order (depending on parameters $d$ and $c$, see \cite{M18}). 
In the present paper, we study a combination of a sublinear weighted function with the max choice from the sample of the growing size. We will show that both a sublinear and linear maximum degree is possible in this case.

Let describe our model. Let fix $\alpha\in (0,1),$  $\gamma\in (0,1)$ and $c_d>0$. Let define $d_n=c_dn^{\gamma}$, $n\in\N$. Let consider i.i.d. random variables $m,\{m_n\}_{n\in\N}$ with values in $N$, such that $\E m^2<\infty$. We would consider sequence of random graph $G_n$, $n\in\Z_+$, that builds by following inductive rule. We start with initial graph $G_1$ that consist of two vertices $v_0$ and $v_1$ and $m_1$ edges between them. Then on $n+1$-th step we add new vertex $v_{n+1}$ and draw $m_{n+1}$ edges from it to vertices $Y_n^1, ...,Y_n^{m_{n+1}}$ choosen from $V(G_n)$ by following rule. For each $i\in N$, $1\leq i \leq m_{n+1}$, we independently (given $G_n$) choose vertices $X_{n}^{i,1},...,X_{n}^{i,d_{n+1}}$ from $V(G_n)$ with probabilities proportional to their degree in power $\alpha$:
$$\Pr(X_{n}^{i}=v_j)=\frac{(\deg_{G_n}v_j)^{\alpha}}{\sum_{k=0}^{n}(\deg_{G_n}v_k)^{\alpha}}.$$
 Then $Y_n^i$ would be the vertex with highest degree among $X_n^{i,1},...,X_n^{i,d_n}$, in case of tie we choose vertex in accordance with fair coin toss.

We would be interested in the number of vertices of fixed degree and the maximum degree of the graph. Let $N_n(k)$ be the number of vertices of degree $k$ in graph $G_n$ and $M_n$ be the maximum degree of vertices in $G_n$. Then the total weight $D_n$ of all vertices in $G_n$ is
$$D_n:=\sum_{i=0}^{n}\left(\deg_{G_n} v_i \right)^{\alpha}=\sum_{k=1}^{\infty}N_{n}(k)k^{\alpha}.$$
There are two ways to increase the maximum degree of the graph. First, we could add new vertex with degree higher then degree of already existing vertices. To prevent that, we would put certain conditions on the tails of $m_n$ which would provide that with high probability $m_n\leq M(n)$ for all large enought $n$. Second, we could increase the maximum degree by drawing edges to the vertex with maximum degree. Given $G_n$, the probability to draw an edge to the vertex with maximum degree equals to
$$\Pr\left(\deg_{G_n}Y_n^i=M(n)\right)=\left(1-\left(1-\frac{(M(n))^{\alpha}L(n)}{D_n}\right)^{d_n}\right).$$
Therefore, evolution of $M(n)$ satisfy
\begin{equation}
\label{eq:M_n}
\begin{gathered}
   \E(M(n+1)-M(n)|\F_n)\geq 
\E m\left(1-\left(1-\frac{(M(n))^{\alpha}}{D_n}\right)^{d_n}\right),\\
\E(M(n+1)-M(n)|\F_n)\leq 
\E m\left(1-\left(1-\frac{(M(n))^{\alpha}L(n)}{D_n}\right)^{d_n}\right),
\end{gathered}
\end{equation}
where $L(n)$ is the number of vertices of degree $M(n)$ and $\F_n$ is sigma-algebra that corresponds to $G_n$.
Let formulate our main results.

\begin{theorem}
\label{th:D_n}
Let $\Pr(m>c)>0$ for any $c>0$. Then
$$\frac{N_k(n)}{n}\to \Pr(m=k)\;a.s.$$
In particular
\begin{equation}
\label{eq:D_n}
    \frac{D(n)}{n}\to \sum_{k=1}^{\infty}k^{\alpha}\Pr(m=k)=\E m^{\alpha}\;a.s.
\end{equation}
\end{theorem}

\begin{theorem}
\label{th:max_degree}
Let $\Pr (m=k)\leq ck^{-\beta}$ for some $\beta>1+\frac{1-\alpha}{\gamma}$ and constant $c>0$. Then
\begin{enumerate}
\item If $\alpha+\gamma<1$, then $\frac{M(n)}{n^{\frac{\gamma}{1-\alpha}}}\to x^{\ast}$, where $x^{\ast}=\left(\frac{\E m (1-\alpha)d_n}{\gamma\E m^{\alpha}}\right)^{\frac{1}{1-\alpha}}$ a.s.
\item If $\alpha+\gamma=1$, then $\frac{M(n)}{n}\to \rho^{\ast}$ a.s., where $\rho^{\ast}$ is a unique positive solution of the equation $1-e^{\frac{c_dx^{\alpha}}{\E m^{\alpha}}}-x$
\item If $\alpha+\gamma>1$, then $\frac{M(n)}{n}\to \E m$ a.s.
\end{enumerate}

\end{theorem}

Theorem~\ref{th:D_n} shows that degrees of most vertices do not change after their appearance. It happens, as would be proven in section~2, due to the increasing size of the sample, which results in vertices with a relatively high degree to be present in the sample with high probability. In other words, the new vertex with high probability connects to the vertices whos degree exceeds a certain growing level. Theorem~\ref{th:max_degree} shows how new edges could be accumulated among vertices with high degrees. In case $\alpha+\gamma>1$ almost all edges would be drawn toward single vertex with degree asymptotically equals to $(\E m) n$, while in the case $\alpha+\gamma<1$ edges would be drawn to the vertices of degrees up to $x^{\ast}n^{\frac{\gamma}{1-\alpha}}$. If we consider $m$ to have power-law distribution then such a combination of max choice with sublinear weighted function would result in vertices with high degrees to follow different exponent then vertices with relatively small degrees up to existing of the condensation for $\alpha+\gamma\geq 1$. 

We will use stochastic approximation techniques to prove almost sure convergence in the linear case. Note that stochastic approximation is widely used to prove almost sure convergence for linear order of maximal degree (see, for example, \cite{MP14,HJ16,HJY18}), while to prove sublinear order of the maximum degree martingale approach is usually used (see, for example, \cite{Mori05,M18}). Also, in contrast with some previous works on models with choice (in particular, \cite{MP14,M18}), due to nonconvexity of the weighted function, we do not use persistent hub argument and instead use auxiliary stochastic approximation processes to separately get lower and upper bound for the maximum degree.

Let give a short description of stochastic approximation approach (for more details see, for example, \cite{Chen03,Pemantle}) that we use to prove our results. Process $Z(n)$ is stochastic approximation process if it could be written as
$$Z(n+1)-Z(n)=\gamma_n\left(F(Z_n)+E_n+R_n\right)$$
where $\gamma_n$, $E_n$ and $R_n$ satisfy following condition. $\gamma_n$ is not random and $\sum_{n=1}^{\infty}\gamma_n>\infty$, $\sum_{n=1}^{\infty}(\gamma_n)^2<\infty$, usially one put $\gamma_n=\frac{1}{n}$ or $\gamma_{n}=\frac{1}{n+1}$. The function $F(x)$ is continues with isolated roots and represent dependence of the increment of the process from its current state. Often process $Z(n)$ belongs to some interval $[a,b]$, and therefore function is considered only on this interval as well. For example, if $Z(n)$ is fraction of balls in urn model it belongs to $[0,1]$. The term $E_n$ is $\F_n$-measurable where $\F_n$ is natural filtration of $Z(n)$, $\E(E_n|\F_n)=0$ and $\E ((E_n)^2|\F_n)<c$ for some fixed constant $c$. Usially one put $E_n=\frac{1}{\gamma_n}(Z(n+1)-\E(Z(n+1)|\F_n))$ and therefore function $F(x)$ could be found from representation $\E(Z(n+1)-Z(n)|\F_n)=\gamma_n(F(Z(n))+R_n)$ where $R_n$ is small error term that satisfy $\sum_{n=1}^{\infty}\gamma_n|R_n|<\infty$ almost surely. If necessary conditions are met, the process will almost surely converge to the zero set of $F(x)$. Moreover, process $Z(n)$ could converge only to stable zero ($x^{\ast}$ is stable zero if $F(x)$ change sigh from $+$ to $-$ when approaching it).  

Note that conditions on $F(x)$ and $R_n$ are usually true and easy to check, while for some representations condition on $E_n$ could break due to multiplication on term $\frac{1}{\gamma_n}$ that turns to infinity. In particular, this is why classical stochastic approximation results could not be used for the sublinear case and instead, we would use a different approach, including large deviation estimates. 
Let us give some outline of this approach. If we consider the degree of certain vertex or the maximal degree of the graph, on step $n$ its increase could be represented as the sum of $m_n$ conditionally independent Bernoulli random variables. Therefore, under certain conditions, we could estimate evolution of the degree by the sum of independent (given the condition) Bernoulli random variables. Then we could consider their expectations and use large deviation results to ensure that the process does to deviate far from its expectation. We would use following standart large deviasion result on bernoulli random variables
\begin{lemma}
\label{lem:large_dev}
Let $\eta_1,...,\eta_n$ be i.i.d. bernoulli variables with parameter $p$. Let $S_n=\sum_{i=1}^n\eta_i$.
Then for any $\delta>0$ there are constants $C$ and $c=c(\delta)>0$, such that for all $n\in\N$ and any $p\in(0,1)$
$$\Pr(S_n\leq (1+\delta) pn)\leq Ce^{-cpn}.$$
\end{lemma}
It's proof well known and uses standard combinatorial argument, so we will not provide it here.

\subsection*{Proof approach and organization}

In section~\ref{sec:number_of_vertices} we prove strong law of large numbers for the number of vertices of fixed and almost sure converges for the total weight of the graph. We would later use it to simplify formulas for stochastic approximation argument.

In section~\ref{sec:lower_bound} we prove Theorem~\ref{th:max_degree} in case $\alpha+\gamma\geq 1$. To do so, we let $L(n)=1$ to approximate the maximum degree from below using stochastic approximation processes. As a result, we would get a linear lower bound for the maximum degree. Then, due to the total degree of the graph being linear, the number of vertices with degrees above linear level is bounded by a constant and hence a simple argument provides that $L(n)=1$ with high probability and therefore we would get almost sure convergence for the case $\alpha+\gamma\geq 1$.

In section~\ref{sec:upper_bound} we provide the proof of Theorem~\ref{th:max_degree} for the case $\alpha+\gamma<1$. We first put  $L(n)=1$ to get lower bound for maximum degree of the graph.
Then we would use a large deviation approach towards possible rate of growth of a fixed vertex to show that with high probability degrees of all vertices do not grow faster than the given rate.

\section{The number of vertices of fixed degree}
\label{sec:number_of_vertices}

In this section, we provide proof of Theorem~\ref{th:D_n}.

\begin{proof}
Note that, since $\Pr(m>c)>0$ for any $c>0$, the number of vertices with degree more then $c$ with high probability of order $n$ for any $c$. Therefore for any $k\in\N$ there is a constant $C_k>0$, such that $D_n-D_{n}(k)\geq C_k n$ with high probability. Since with high probability $D_n\leq 2n\E m $, we get that there is a constant $c_k\in(0,1)$, such that with high probability $\frac{D_{n}(k)}{D_n}\leq c_k$. Hence, with high probability
\begin{align*}
\E\left(\1\{\deg Y_{n}^i=k\}|\F_n\right)
&=\left(\sum_{j=1}^{k}\frac{N_{n}(j)j^{\alpha}}{D_n}\right)^{d_n} -\left(\sum_{j=1}^{k-1}\frac{N_{n}(j)j^{\alpha}}{D_n}\right)^{d_n}\\
&=\exp\left\{d_{n}\ln\left(\frac{D_n(k)}{D_n}\right)\right\} -\exp\left\{d_n\ln\left(\frac{D_n(k-1)}{D_n}\right)\right\}\\
&\leq \exp\{d_n\ln(c_k)\}\to 0
\end{align*}
as $n\to\infty$.
Therefore, almost all vertices with degree $k$ do not have edges drawn into them which results in the first statement of the theorem. To get the second statement note that the total weight of vertices with degrees more then $k$ at time $n$ does not exceed $\frac{\sum_{i=1}^{n}m_i}{k}k^{\alpha}=n\E(m)o_k(1)$ a.s. as $k\to\infty$.
Hence
$$\frac{D_n}{n}=\sum_{j=1}^{k}\frac{k^{\alpha}N_k(n)}{n}+\E mo_k(1)\to\sum_{j=1}^{\infty}\frac{k^{\alpha}N_k(n)}{n}$$
a.s. as $k\to\infty$.
Since $\frac{N_k(n)}{n}\to \Pr(m=k)$ a.s. as $n\to\infty$, we get that
$$\frac{D_n}{n}\to\sum_{j=1}^{\infty}k^{\alpha}\Pr(m=k)=\E m^{\alpha}$$
a.s. as $n\to\infty$.
\end{proof}

\section{The maximum degree: case $\alpha+\gamma\geq 1$}
\label{sec:lower_bound}

First, we provide an estimate of the maximum degree from below.
To do so we would put $L(n)=1$ in formula \eqref{eq:M_n}. Then we would use a stochastic approximation to prove the convergence of the resulting process.

We get that
\begin{align*}
\frac{1}{\E m }\E(M(n+1)-M(n)|\F_n)&\geq \left(1-\left(1-\frac{(M(n))^{\alpha}}{D_n}\right)^{d_n}\right)\\
&=\left(1-\exp\left\{d_n\ln\left(1-\frac{(M(n))^{\alpha}}{D_n}\right)\right\}\right)\\
&\geq \left(1-\exp\left\{-\frac{d_n(M(n))^{\alpha}}{D_n}\right\}\right)\\
&=1-\exp\left\{-\frac{c_d\left(\frac{M(n)}{n}\right)^{\alpha}}{\E m^{\alpha}n^{1-\gamma-\alpha}}(1+o(1))\right\} \\
&\geq 1-\exp\left\{-\frac{c_d}{\E m^{\alpha}}\left(\frac{M(n)}{n}\right)^{\alpha}(1+o(1))\right\}.
\end{align*}
There exists $A(n)=A(n,\epsilon,n_0)$ such that
$A(n_0)=M(n_0)$, 
$$\E\left(A_{\epsilon}(n+1)-A_{\epsilon}(n)|F_n\right):=\E m \left(1-\exp\left\{-\frac{c_d}{\E m^{\alpha}}\left(\frac{A(n)}{n}\right)^{\alpha}\right\}\right)$$
and $A(n)\leq M(n)$ on $\mathcal{A}_{\epsilon}(n_0)$. Consider $B(n):=A(n)/n$.
Then
$$\E(B(n+1)-B(n)|\F_n)=\frac{\E m }{n+1}\left(1-\exp\left\{-\frac{c_d}{\E m^{\alpha}}\left(B(n)\right)^{\alpha}\right\}-B(n)\right).$$
Note that function $g(x)=1-e^{-\frac{c_dx^{\alpha}}{\E m^{\alpha}}}-x$ has a unique and stable root $\rho^{\ast}$ in $[0,1]$. Also, $|n(B(n+1)-B(n))|\leq m_{n+1}$ and hence 
$$\E((n(B(n+1)-B(n)))^2|\F_n)\leq \E m^2<\infty.$$ 
Therefore
$B(n)\to \rho^{\ast}$ a.s. as $n\to\infty$.
As result, for $\alpha+\gamma\geq 1$ we get that $\liminf \frac{M(n)}{n}\geq \rho^{\ast}$ a.s. Moreover, for $\alpha+\gamma> 1$ we get that
$$\E(M(n+1)-M(n)|\F_n)\geq \E m \left(1-\exp\left\{-\frac{c_d\left(\rho^{\ast}\right)^{\alpha}}{\E m^{\alpha}}(1+o(1))n^{\gamma+\alpha-1}\right\}\right)\to \E m $$
a.s., and therefore $\frac{M(n)}{n}\to \E m$ a.s.
Note that in the case $\alpha+\gamma=1$ putting $L(n)=1$ would give us actual bound. Indeed, since $\liminf_{n\to\infty}\frac{M(n)}{n}=\rho^{\ast}>0$ (which corresponds to $\alpha+\gamma\geq 1$), we get that with high probability $L(n)\leq \frac{\E m}{\rho^{\ast}}$. It is a well-known fact that the probability for the simple random walk to return to the origin is $O(n^{-1/2})$. Adding the positive probability to not move does not change that asymptotic. Note that the probability to increase the degree of the vertex with a degree higher the $\rho^{\ast}n$ is bound from below by some constant. Also, for two vertices with different degrees probability to increase degree is higher for the vertex with a higher degree. As a result, for any pair of vertices with degrees higher then $\rho^{\ast}n$ the probability that they have the same degree is at time $n$ is $O(n^{-1/2})$.  Hence $\Pr(L(n)>1)=O(n^{-1/2})$ and therefore $\frac{M(n)}{n}\to\rho^{\ast}$ a.s. 

\section{The maximum degree: case $\alpha+\gamma< 1$}
\label{sec:upper_bound}
   
In this section, we prove Theorem~\ref{th:max_degree} in the case $\alpha+\gamma< 1$. First, similar to the case $\alpha+\gamma\geq 1$, we put $L(n)=1$ in formula \eqref{eq:M_n} to get lower bound. To get a matching upper bound we consider the evolution of the degree of a fixed vertex and get a large deviation type estimate for it. We would prove that probability to grow faster than a certain rate has an exponential tail and therefore with high probability no vertices degrees grow faster than this rate.
Once again, we get that
\begin{align*}
\E(M(n+1)-M(n)|\F_n)&\geq \E m\left(1-\exp\left\{-\frac{d_n(M(n))^{\alpha}}{D_n}\right\}\right)\\
&= \E m\left(\frac{d_n(M(n))^{\alpha}}{D_n}+O\left(\left(\frac{d_n(M(n))^{\alpha}}{D_n}\right)^2\right)\right).
\end{align*}
Note that, due to Theorem~\ref{th:D_n}, for any $\epsilon>0$ probability of event $\mathcal{A}_{\epsilon}(n)=\{\forall l\geq n: D_l<\E m^{\alpha}+\epsilon\}$ turns to $1$ as $n$ turns to $\infty$. Recall that $d_n=c_dn^{\gamma}$.
Therefore for any $n_0\in\N$, if $\frac{d_n}{n^{1-\alpha}}=c_dn^{\gamma+\alpha-1}\to 0$ (meaning $\gamma+\alpha<1$) then it is possible to define a process $A_{\epsilon}(n)$, $n\leq n_0$, such that $A_{\epsilon}(n_0)=M(n_0)$,
$$\E\left(A_{\epsilon}(n+1)-A_{\epsilon}(n)|F_n\right):=\E m \frac{d_n(A_{\epsilon}(n))^{\alpha}}{(\E m^{\alpha}+\epsilon)n}=\frac{c_d\E m(A_{\epsilon}(n))^{\alpha}}{(\E m^{\alpha}+\epsilon)n^{1-\gamma}}$$
and $A_{\epsilon}(n)\leq M(n)$ on $\mathcal{A}_{\epsilon}(n_0)$.
Let $x^{\ast}:=\left(\frac{c_d\E m (1-\alpha)}{\gamma\E m^{\alpha}}\right)^{\frac{1}{1-\alpha}}$. Fix $\delta>0$. Consider event $\mathcal{C}_n=\mathcal{C}_n(\epsilon,\delta)=\{A_{\epsilon}(n)\leq (1-\delta)x^{\ast}n^{\frac{\gamma}{1-\alpha}}\}$. On this event
\begin{align*}
\E\left(A_{\epsilon}(n+1)-A_{\epsilon}(n)|F_n\right)&\geq\frac{c_d\E m((1-\delta)x^{\ast}n^{\frac{\gamma}{1-\alpha}})^{\alpha}}{(\E m^{\alpha}+\epsilon)n^{1-\gamma}}\\
&=\frac{(c_d\E m)^{\frac{1}{1-\alpha}}(1-\delta)^{\alpha}(1-\alpha)^{\frac{\alpha}{1-\alpha}}}
{\gamma^{\frac{\alpha}{1-\alpha}}(\E m^{\alpha})^{\frac{\alpha}{1-\alpha}}(\E m^{\alpha}+\epsilon)n^{1-\frac{\gamma}{1-\alpha}}}\\
&=\frac{\gamma(1-\delta)^{\alpha}x^{\ast}}{(1-\alpha)\left(1+\frac{\epsilon}{\E m^{\alpha}}\right)n^{1-\frac{\gamma}{1-\alpha}}}.
\end{align*}
Note that
$$\sum_{n=n_1}^{n_2}\frac{\gamma(1-\delta)^{\alpha}x^{\ast}}{(1-\alpha)\left(1+\frac{\epsilon}{\E m^{\alpha}}\right)n^{1-\frac{\gamma}{1-\alpha}}}
\geq\frac{(1-\delta)^{\alpha}x^{\ast}}{\left(1+\frac{\epsilon}{\E m^{\alpha}}\right)}\left(n_2^{\frac{\gamma}{1-\alpha}}-n_1^{\frac{\gamma}{1-\alpha}}\right).$$
Choose $\delta$ and $\epsilon$ such that $\frac{(1-\delta)^{\alpha}}{\left(1+\frac{\epsilon}{\E m^{\alpha}}\right)}>(1-\delta)$. Let put $\sigma:=\frac{(1-\delta)^{\alpha}}{\left(1+\frac{\epsilon}{\E m^{\alpha}}\right)}-(1-\delta)>0$.
Note that $A_{\epsilon}(n+1)-A_{\epsilon}(n)$ is the sum of $m_{n+1}$ Bernoulli random variables. Therefore $A_{\epsilon}(n_2+1)-A_{\epsilon}(n_1)$ is the sum of $\sum_{n=n_1+1}^{n_2}m_n$ Bernoulli random variables and hence, on $\cup_{n_1\leq n\leq n_2}\mathcal{C}_n$, by large deviation estimate we get that
\begin{align*}
\Pr\left(A_{\epsilon}(n_2+1)-A_{\epsilon}(n_1)<(1-\delta+\sigma/2)x^{\ast}\left(n_2^{\frac{\gamma}{1-\alpha}}-n_1^{\frac{\gamma}{1-\alpha}}\right)\right)&\leq Ce^{-cn^{\gamma-1}\sum_{n=n_1+1}^{n_2}m_n}\\
&\leq Ce^{-cn^{\gamma}}
\end{align*}
for some constants $c,C>0$. Hence process $A(n)$ with high probability could not stay below level $(1-\delta)x^{\ast}n^{\frac{\gamma}{1-\alpha}}$ which gives us estimate
$$\liminf_{n\to\infty}M(n)\geq\liminf_{n\to\infty}A_{\epsilon}(n)\geq (1-\delta)x^{\ast}n^{\frac{\gamma}{1-\alpha}}$$
almost surely, and therefore
$$\liminf_{n\to\infty}M(n)\geq x^{\ast}n^{\frac{\gamma}{1-\alpha}}$$
almost surely.

Now, let prove matching upper bound.

Note that for any $n_0>0$ and $\epsilon\in (0,\frac{\gamma}{1-\alpha})$ 
\begin{align*}
\Pr(\exists n\geq n_0: m_{n}\geq n^{\frac{\gamma}{1-\alpha}-\epsilon})&\leq \sum_{n=n_0}^{\infty}\sum_{i=n^{\frac{\gamma}{1-\alpha}-\epsilon}}^{\infty}\pr(m_n=i)\\
&=\sum_{n=n_0}^{\infty}cn^{(1-\beta)\left(\frac{\gamma}{1-\alpha}-\epsilon\right)}\to 0
\end{align*}
as $n_0\to\infty$ if $(1-\beta)\left(\frac{\gamma}{1-\alpha}-\epsilon\right)<-1$. Note that such $\epsilon$ exists since $\beta>1+\frac{\gamma}{1-\alpha}$.

Let consider events $\C_{n_0}=\{\exists \epsilon>0:\forall n\geq n_0\; m_n\leq n^{\frac{\gamma}{1-\alpha}-\epsilon} \}.$ Then it is enought to prove the result on $\C_{n_0}$ for any fixed $n_0$.

On $\C_{n_0}$ we get that
$$\E(M(n+1)-M(n)|\F_n)= \E m\left(1-\left(1-\frac{L_n(M(n))^{\alpha}}{D_n}\right)^{d_n}\right).$$

Let estimate from above the condition probability $p_n(v)$ to draw an edge to a single vertex. Note that such probability is increasing under the condition that there are no vertices with a higher degree in the sample, which achieved for vertices with the highest degree. For condition probability $q_n$ to draw an edge to a certain vertex with the highest degree we get
\begin{align*}
q_n&=\frac{1}{L_n}\left(1-\left(1-\frac{L_n(M(n))^{\alpha}}{D_n}\right)^{d_n}\right)\\
&\leq\frac{1}{L_n}\frac{L_n(M(n))^{\alpha}d_n}{D_n}=\frac{(M(n))^{\alpha}d_n}{D_n}.
\end{align*}
Therefore for any vertices we have
$$p_n(v)\leq \frac{(\deg_n(v))^{\alpha}d_n}{D_n}.$$
Let fix $\epsilon>0$. Introduce event $\B=\B(n_0,\epsilon)=\{\forall n\geq n_0:D_n\leq \frac{(\E m^{\alpha})n}{1+\epsilon}\}$. Then $\Pr(\B)\to 1$ as $n_0\to\infty$. Then on $\B$
$$p_n(v)\leq (1+\epsilon)\frac{(\deg_n(v))^{\alpha}}{(\E m^{\alpha})n^{1-\gamma}}.$$
Therefore for any $n_1\geq n_0$ the evolution of degree $\deg_n(v)$ for $n\geq n_1$ could be estimated from above by the sum of $m_n$ bernoulli random variables $\xi_{i,n},$ $i=1,...,m_n$ that build using i.i.d. random variables $u_{i,n},$ $i\in\N,n\in\N$ that uniformal on $[0,1]$ as follow $$\xi_{i,n}=\1\left\{u_{i,n}\leq(1+2\epsilon)\frac{(S_n)^{\alpha}}{(\E m^{\alpha})n^{1-\gamma}}\right\},$$
where $S_n=\deg_{n_1}(v)+\sum_{k=n_1}^{n-1}\sum_{i=1}^{m_k}\xi_{i,k}.$ Let define $\tau_n:=\sum_{i=1}^{n}m_i$ and $\pi_n:=\inf{m:\tau_n=m}$.
Since probability in the right side is increasing when $S_n$, if we instead making $m_n$ steps at one moment consider $m_n$ consequtive steps we would increase corresponding probabilities. Hence, the evolution of degree $\deg_n(v)$ for $n\geq n_1$ could be dominated by $Y_{\tau_n}$ where $Y_{n+1}-Y_n$ are bernoulli random variables that satisfy
$$Y_{n+1}-Y_n=\1\left\{u_n\leq(1+2\epsilon)\frac{(Y_n)^{\alpha}}{(\E m^{\alpha})(\pi_n)^{1-\gamma}}\right\}$$
for some i.i.d uniformal on $[0,1]$ random variables $u_n$, $n\in\N$.
Note that since $\beta>2$ we get that $\frac{\tau_{n}}{n}\to E_m$ a.s. Therefore for any $\delta>0$ 
$$\Pr \left(\left|\frac{\tau_{n}}{n}-\E m\right|>\delta\right)\to 0$$
as $n\to\infty$. Hence if we consider event $\mathcal{E}_{n_0}=\left\{\forall n>n_0:(\pi_n)^{1-\gamma}<\frac{2+\epsilon}{3+\epsilon}( n/\E m)^{1-\gamma}\right\}$ we get that $\Pr\left(\mathcal{E}_{n_0}\right)\to 0$
as $n_0\to\infty$.
On $\mathcal{E}_{n_0}$ for $n\geq n_0$ variables $Y_n$ could be dominated by bernoulli random variables $X_n$, such that
$$\zeta_n:=X_{n+1}-X_n=\1\left\{u_n\leq(1+3\epsilon)\frac{(X_n)^{\alpha}}{(\E m^{\alpha})(n/\E m)^{1-\gamma}}\right\}.$$

Using lemma~\ref{lem:large_dev} we get an estimate on the growth rate of $X_n$. Fix $\sigma>0$. For probablity that $X_{n(1+\delta)}\geq (1+\delta)^{\frac{\gamma}{1-\alpha}}X_{n}$ on event $\{X_n\geq (1+\sigma)(1+3\epsilon)\left(\frac{(\E m)^{1-\gamma} (1-\alpha)}{\gamma\E m^{\alpha}}\right)^{\frac{1}{1-\alpha}}n^{\frac{\gamma}{1-\alpha}}\}$ we get
\begin{align*}
\Pr\!\left(X_{n(1+\delta)}\!\geq\! (1+\delta)^{\frac{\gamma}{1-\alpha}}X_{n}\right)\! &=\!\Pr\left(\sum_{i=n}^{(1+\delta)n}\zeta_i\geq\left((1+\delta)^{\frac{\gamma}{1-\alpha}}-1\right)X_n\right)\\
&=\!\Pr\!\left(\!\sum_{i=n}^{(1+\delta)n}\!(\zeta_i\!-\!\E\zeta_i) \!\geq\!\left((1+\delta)^{\frac{\gamma}{1-\alpha}}\!-\!1\right)X_n\!-\!\sum_{i=n}^{(1+\delta)n}\!\E\zeta_i\!\right)
\end{align*}
Note that 
$$(1+\delta)^{\frac{\gamma}{1-\alpha}}X_{n}^{1-\alpha}-\frac{\delta n(1+3\epsilon)}{(\E m^{\alpha})(n/\E m)^{1-\gamma}}\geq$$
$$\geq\left((1+\delta)^{\frac{\gamma}{1-\alpha}}-1\right) \left((1+\sigma)(1+3\epsilon)\left(\frac{(\E m)^{1-\gamma} (1-\alpha)}{\gamma\E m^{\alpha}}\right)^{\frac{1}{1-\alpha}}n^{\frac{\gamma}{1-\alpha}}\right)^{1-\alpha}$$
$$-\frac{\delta (1+3\epsilon)(\E m)^{1-\gamma}}{(\E m^{\alpha})}n^{\gamma}$$
$$=\left(\left((1+\delta)^{\frac{\gamma}{1-\alpha}}-1\right)(1+\sigma)^{1-\alpha}\frac{1-\alpha}{\gamma}(1+3\epsilon)^{-\alpha}-\delta \right)\left(\frac{(1+3\epsilon)(\E m)^{1-\gamma}}{(\E m^{\alpha})}\right)n^{\gamma}
$$
$$>c\frac{\delta n(1+3\epsilon)}{(\E m^{\alpha})(n/\E m)^{1-\gamma}}$$
for some $c=c(\delta,\epsilon)>0$ for small enought $\delta$ and $\epsilon$. Hence 
\begin{align*}
\Pr\left(X_{n(1+\delta)}\geq (1+\delta)^{\frac{\gamma}{1-\alpha}}X_{n}\right) &\leq\Pr\left(\sum_{i=n}^{(1+\delta)n}(\zeta_i-\E\zeta_i) 
\geq c\sum_{i=n}^{(1+\delta)n}\E\zeta_i\right)\\
&\leq Ce^{-c_1n^{\gamma}}
\end{align*}
for some $C,c_1>0$ that does not depend on $n$. Therefore probability that at any time $n\geq n_0$ the vertex with degree $w(n)\geq(1+\sigma)(1+3\epsilon)\left(\frac{(\E m)^{1-\gamma} (1-\alpha)}{\gamma\E m^{\alpha}}\right)^{\frac{1}{1-\alpha}}n^{\frac{\gamma}{1-\alpha}}$ would have degree more then $(1+\delta)^{\frac{\gamma}{1-\alpha}}w(n)$ at time $(1+\delta)n$ turns to $0$ as $n_0$ turns to $\infty$.
Therefore
$$\Pr\!\left(\!\forall n\geq n_0\,:\,M(n)\geq(1\!+\!\delta)^{\frac{\gamma}{1-\alpha}}(1\!+\!\sigma)(1\!+\!3\epsilon)\left(\frac{(\E m)^{1-\gamma} (1-\alpha)}{\gamma\E m^{\alpha}}\right)^{\frac{1}{1-\alpha}}n^{\frac{\gamma}{1-\alpha}}\!\right)\!\to 0$$ as $n_0\to\infty$.
As result we get that for any $\sigma>0$ we could find small enought $\delta$ and $\epsilon$, such that
$$\limsup_{n\to\infty}\frac{M(n)}{n^{\frac{\gamma}{1-\alpha}}}\leq (1+\delta)^{\frac{\gamma}{1-\alpha}}(1+\sigma)(1+3\epsilon)\left(\frac{(\E m)^{1-\gamma} (1-\alpha)}{\gamma\E m^{\alpha}}\right)^{\frac{1}{1-\alpha}}.$$
Therefore
$$\limsup_{n\to\infty}\frac{M(n)}{n^{\frac{\gamma}{1-\alpha}}}\leq\left(\frac{(\E m)^{1-\gamma} (1-\alpha)}{\gamma\E m^{\alpha}}\right)^{\frac{1}{1-\alpha}},$$
which concludes the proof of the upper bound for the case $\alpha+\gamma<1$.

\section{Discussion}
In present paper we achived transition in behavior of the maximum degree between sublinear case of $M(n)\sim x^{\ast}n^{\frac{\gamma}{1-\alpha}}$ and strict linear case $M(n)\sim (\E m) n$ with transition of the type $M(n)\sim \rho^{\ast}n$. We considered added $m$ edges on each step with tails of $m$ bounded from above by power law with power more then $1+\frac{1-\alpha}{\gamma}$, and therefore tails of $m$ did not affect th maximum degree. One could consider case of power law tails of $m$ with power $\beta<1+\frac{1-\alpha}{\gamma}$, it seems that argument given in prove of the upper bound sshould provide that once vertex with degree above level $x^{\ast}n^{\frac{\gamma}{1-\alpha}}$ emerges, it could not keep up with this level and its degree should turn to $n^{\frac{\gamma}{1-\alpha}}$ as $n\to\infty$. 

Note that considering $d_n=c_dn^{\gamma}$ gives us maximum degree of order at least $n^{\gamma}$ even without considering preferential attachment (if we put $\alpha=0$). It could be interesting to see if transition between $(\ln n)^{b}$ and $n^{a}$ orders (for $b>0$, $0<a<1$) of the could be found if we consider $d_n$ of order $(\ln n)^{c_1}$ and weighted function of the type $\frac{x}{(\ln x)^{c_2}}$, $c_1,c_2>0$. 

The other modification of the model is to consider the combination of the min choice with superlinear function. It is not clear if the power-law type of maximum degree could be achieved in this case. For example, in \cite{HJ16} for meek choice (when we choose vertex with $s-th$ highest degree for $s>1$) was shown that the maximum degree could be either of linear order or of $(\ln n)^b$ order with no power-law type of behavior.

\section*{Acknowledgements.}
The presented work was funded by a grant from the Russian Science Foundation (project No. 19-71-00043).

\bibliographystyle{alpha}
\bibliography{nonlinear preferential attachment with choice}

\end{document}